# DENSITY OF CUBIC FIELD DISCRIMINANTS

DAVID P. ROBERTS

ABSTRACT. In this paper we give a conjectural refinement of the Davenport-Heilbronn theorem on the density of cubic field discriminants. We explain how this refinement is plausible theoretically and agrees very well with computational data.

1. INTRODUCTION

Let $a_n$ be the number of isomorphism classes of abelian cubic fields with discriminant $n$. Let $b_n$ be the number of isomorphism classes of non-abelian cubic fields with discriminant $n$. The numbers $a_n$ are very well understood. The numbers $b_n$ have been the subject of extensive theoretical and computational study for at least sixty years, but are less well understood. The object of this note is to contribute to the study of these $b_n$, by bringing together the theoretical and computational literature.

For $\alpha \in \{-, +\}$, define

$$(1.1) \qquad g_\alpha(x) = \sum_{n=1}^{x} b_{\alpha n}.$$

Put $C_- = 3/4$ and $C_+ = 1/4$. The main theoretical result is due to Davenport and Heilbronn [8]:

$$(1.2) \qquad g_\alpha(x) \sim C_\alpha \frac{1}{3\zeta(3)} x.$$

The numbers $b_{\alpha n}$ have been computed for larger and larger $n$ by many authors. The most recent results are due to Belabas [3], who introduced a vastly improved method and worked through the cutoff $x = 10^{11}$.

It has been discussed by several authors that the "match" between computational evaluation of $g_\alpha(x)$ and the asymptotic formula (1.2) is not very good, as illustrated by Tables 1, 2, and 3 below. For example, the paper [9] dealt with the case $\alpha = -$; here is an excerpt from pages 322-324, with some trivial notational changes:

> If, however, the reader were to plot the density versus $x$, he would be somewhat astonished to see that this density is increasing so slowly that his first impression would be that it will not make it to Davenport-Heilbronn limit. Thus it remains a challenging problem, assuming that the D-H limit is not in error, to explain the origin of this very slow convergence. This problem was indicated by Shanks in [13],[1], and on the real side in [13],[10], where the problem is further aggravated by even slower convergence. To date, and to our knowledge, no good quantitative explanation of this phenomenon has been given.

In this paper we conjecture a quantitative explanation of this slow-convergence phenomenon.

---

## 2. The abelian case

It is best to briefly review the situation in the abelian case. One has the formula

$$(2.1) \qquad \sum_{n=1}^{\infty} \frac{a_n}{n^s} = -1 + \left(1 + \frac{2}{81^s}\right) \prod_{p \equiv 1\ (6)} \left(1 + \frac{2}{p^{2s}}\right).$$

Put

$$(2.2) \qquad f(x) = \sum_{n=1}^{x} a_n.$$

Cohn [4] showed that the function (2.1) continues analytically to $\mathrm{Re}(s) > 1/4$, except for a simple pole at $s = 1/2$. He computed the residue of this pole. This gave him more than enough information to conclude that

$$(2.3) \qquad f(x) \sim \left(\frac{11\sqrt{3}}{36\pi} \prod_{p \equiv 1\ (6)} \frac{(p+2)(p-1)}{p(p+1)}\right) x^{1/2}.$$

Define $h_-(x) = g_-(x)$ and

$$(2.4) \qquad h_+(x) = g_+(x) + \frac{1}{3} f(x).$$

From (1.2) and (2.3) one sees that (1.2) would remain true with $g_\alpha$ replaced by $h_\alpha$. We henceforth focus more on $h_\alpha$, as it appears naturally in (3.2) and (3.3).

## 3. The basic conjecture

Belabas proved [3, (8) and (20)] that, for either $\alpha$, the difference between the sides of (1.2) is $o(x/(\log x)^2)$. From his extensive computational data he speculated [3, page 1232] that these differences are $o(x/(\log x)^\beta)$ for arbitrarily large $\beta$. We believe a stronger statement. Put $K_- = 3/(3+\sqrt{3})$ and $K_+ = \sqrt{3}/(3+\sqrt{3})$.

**Conjecture.** *For each $\alpha$,*

$$(3.1) \qquad h_\alpha(x) = C_\alpha \frac{1}{3\zeta(3)} x + K_\alpha \frac{\zeta(1/3)3(3+\sqrt{3})\Gamma(1/3)^3}{10\pi^3 \zeta(2)\zeta(5/3)} x^{5/6} + o(x^{5/6}).$$

The heuristic we are about to sketch is too crude to suggest any particular error term. The $o(x^{5/6})$ we have written corresponds to the weakest conjecture in which our conjectured new term actually plays a role. However, Tables 1 and 2 below suggest that the error term should be reducible at least most of the way to $o(x^{1/2})$.

We are led to our conjecture theoretically as follows. Define a Dirichlet series

$$(3.2) \qquad \xi_\alpha(s) \;=\; \sum_{k \in \mathcal{K}_\alpha} \frac{1}{|\mathrm{Aut}(k)||D_k|^s}.$$

Here $\mathcal{K}_\alpha$ is the set of isomorphism classes of cubic fields $k$ with discriminant of sign $\alpha$. Also, $D_k$ is the discriminant of $k$ while $\mathrm{Aut}(k)$ is the automorphism group of $k$; so $|\mathrm{Aut}(k)|$ is 3 or 1 according to whether $k$ is abelian or not, the former case occuring only for $\alpha = +$.

For $S$ a finite set of primes, define a weighted version of (3.2)

$$(3.3) \qquad \xi_{\alpha,S}(s) = \sum_{k \in \mathcal{K}_\alpha} \frac{\eta_{k,S}(s)}{|\mathrm{Aut}(k)||D_K|^s}.$$

Here the weights $\eta_{k,S}(s)$ are given on page 122 of [6]. These weights are themselves Dirichlet series $\sum c_m(k,S) m^{-s}$ with integral coefficients $c_m(k,S)$ decreasing to the Kronecker symbol $\delta_{m1}$ as $S$ grows. Thus $\xi_{\alpha,S} \to \xi_\alpha$ in the sense of formal Dirichlet series.

The point of introducing the approximations $\xi_{\alpha,S}$ is that one understands them analytically better than one understands $\xi_\alpha$ itself. From Davenport's earlier work [7] one knows that the rightmost pole of $\xi_{\alpha,\emptyset}$ is at $s = 1$. The analytic behavior of $\xi_{\alpha,\emptyset}$ was investigated more explicitly by Shintani [14], who found the next pole at $s = 5/6$. More generally, it is proved in [5, Prop. 6.2], that each $\xi_{\alpha,S}$ continues to the whole $s$-plane, meromorphically with only poles at $s = 1$ and $s = 5/6$. These poles are simple, their residues $C_\alpha r_{1,S}$ and $K_\alpha r_{5/6,S}$ being deducible from the first two equations in [5, Thm. 6.2i]. In the case $s = 1$, one has to set $b_\alpha$ equal to zero in the corresponding equation; this $b_\alpha$ corresponds to quadratic fields, included in $\xi_{\alpha,S}$ in [5] and [6], but not included here.

From the formulas in [5, Thm. 6.2i], one sees that the numbers $r_{1,S}$ are positive, and decreasing as $S$ grows. Similarly, the numbers $r_{5/6,S}$ are negative and increasing as $S$ grows. Let $r_1$ and $r_{5/6}$ be their respective limits. One does not have enough control at present to rigorously apply any "poles-control-growth" theorems. However the general formalism suggests that $h_\alpha(x)$ grows as $C_\alpha r_1 x + K_\alpha(6/5) r_{5/6} x^{5/6}$. This is the source of our Conjecture 3.1. The constant $r_{5/6}$ is actually nicer than we had initially expected from its source as the limit of the $r_{5/6,S}$; special zeta values arise because of the factorization of the right side of (5.3) below.

To follow our heuristic in a more detailed way, one should first understand it in the more refined context of Section 5 below. We recommend looking next at Prop. 2.1 and the discussion around it in [6]; this treats in more classical language the special case of [5, Thm 6.2i] we need, except it does not repeat the formula for $r_{5/6,S}$.

## 4. Numerical evidence for the basic conjecture

Our conjecture for $\alpha = -$ compares with Belabas's data [3, page 1232] for negative-discriminant cubic fields as in Table 1. Here and below, $x = 10^j$. Also $H_\alpha(x)$ denotes the right side of (1.2) while $H_\alpha^*(x)$ denotes the first two terms of the right side of (3.1). Clearly $H_-^*$ matches the data $h_-$ substantially better than $H_-$ does. The poor match of $H_-$ to $h_-$ in the range $10^5 \leq x \leq 10^6$ formed the basis of the quote in Section 1.

Table 1. Computation and theory: negative-discriminant cubic fields

| $j$ | $h_-(x)$ | $\dfrac{h_-(x)}{H_-(x)}$ | $\dfrac{h_-(x)}{H_-^*(x)}$ | $\dfrac{h_-(x) - H_-^*(x)}{x^{1/2}}$ |
|---|---|---|---|---|
| 2 | 7 | 0.337 | 0.7843510 | -0.1925 |
| 3 | 127 | 0.611 | 0.9993210 | -0.0027 |
| 4 | 1,520 | 0.731 | 0.9943300 | -0.0867 |
| 5 | 17,041 | 0.819 | 0.9998781 | -0.0066 |
| 6 | 182,417 | 0.877 | 1.0001096 | 0.0200 |
| 7 | 1,905,514 | 0.916 | 1.0000100 | 0.0060 |
| 8 | 19,609,185 | 0.943 | 0.9999394 | -0.1188 |
| 9 | 199,884,780 | 0.961 | 0.9999850 | -0.0951 |
| 10 | 2,024,660,098 | 0.974 | 1.0000009 | 0.0176 |
| 11 | 20,422,230,540 | 0.982 | 1.0000003 | 0.0218 |

For $\alpha = +$, we start from the values of $g_+(x) + f(x)$ listed on [3, page 1232]. We use (2.1) to compute $f(x)$ and thereby deduce $g_+(x)$. Our hope that the error term in Conjecture 3.1 can be reduced to $o(x^{1/2})$ would require that the numbers in the last column of Table 1 and the next to last column of Table 2 tend to zero. The data seems rather borderline in support of $o(x^{1/2})$, but strongly indicative that our conjecture as stated is true, and moreover that the $o(x^{5/6})$ can be reduced at least somewhat.

Table 2. Computation and theory: positive-discriminant cubic fields

| $j$ | $g_+(x)$ | $f(x)$ | $\dfrac{h_+(x)}{H_+(x)}$ | $\dfrac{h_+(x)}{H_+^*(x)}$ | $\dfrac{h_+(x)-H_+^*(x)}{x^{1/2}}$ | $\dfrac{g_+(x)-H_+^*(x)}{x^{1/2}}$ |
|---|---|---|---|---|---|---|
| 2 | 0 | 2 | 0.096 | 8.5889786 | 0.0589 | -0.0078 |
| 3 | 22 | 5 | 0.341 | 1.0461129 | 0.0330 | -0.0197 |
| 4 | 366 | 16 | 0.536 | 0.9900166 | -0.0374 | -0.0908 |
| 5 | 4,753 | 51 | 0.688 | 1.0010833 | 0.0163 | -0.0374 |
| 6 | 54,441 | 159 | 0.786 | 0.9988436 | -0.0631 | -0.1161 |
| 7 | 592,421 | 501 | 0.855 | 0.9999134 | -0.0162 | -0.0690 |
| 8 | 6,246,698 | 1,592 | 0.901 | 1.0000259 | 0.0161 | -0.0369 |
| 9 | 64,654,353 | 5,008 | 0.933 | 1.0000097 | 0.0198 | -0.0330 |
| 10 | 661,432,230 | 15,851 | 0.954 | 0.9999988 | -0.0081 | -0.0609 |
| 11 | 6,715,773,873 | 50,152 | 0.969 | 1.0000002 | 0.0046 | -0.0482 |

From Cohn's result (2.3), at most one linear combination $g_+(x) + uf(x)$ can agree with $H_+^*(x)$ to within $o(x^{1/2})$. The last two columns of Table 2 suggest that $u = 1/3$ is more likely than $u = 0$. In fact $u = 0.31\ldots$ would give the best fit, in a least squares sense, and this agrees well with the theoretically favored $u = 1/3$.

## 5. The refined conjecture

The Davenport-Heilbronn theorem comes in a more refined version and so does our conjecture. Let $A$ be the set of symbols $\{111, 21, 3, 1^2 1, 1^3\}$. Then for $K$ a cubic field and $p$ a prime, the splitting behavior of $p$ in $K$ determines an element of $A$. If $p$ is unramified, then it either splits completely (111), splits partially (21), or remains inert (3). If $p$ is ramified, then it is either partially ramified ($1^2 1$) or totally ramified ($1^3$).

Up until now, we have been using $\alpha$ to denote a sign, $+$ or $-$. Henceforth, we use $\alpha$ to denote the following data. First, a finite set $\mathrm{support}(\alpha) \subset \{\infty, 2, 3, 5, \ldots\}$ of places of $\mathbf{Q}$. Second, if $\infty \in \mathrm{support}(\alpha)$, a sign $\alpha_\infty \in \{+, -\}$, as before. Third, for each prime $p \in \mathrm{support}(\alpha)$, an element $\alpha_p \in A$. Let $\mathcal{K}_\alpha$ be the set of isomorphism classes of cubic fields meeting the given local condition $\alpha_v$ for each $v \in \mathrm{support}(\alpha)$. One naturally has $f_\alpha$, $g_\alpha$, $h_\alpha$, $\xi_\alpha$, and $\xi_{\alpha,S}$, all direct generalizations of the notions introduced above in the special case $\mathrm{support}(\alpha) = \{\infty\}$.

Define $C_{\infty,-} = 3/4$ and $C_{\infty,+} = 1/4$, as in Section 1, and $K_{\infty,-} = 3/(3+\sqrt{3})$ and $K_{\infty,+} = \sqrt{3}/(3+\sqrt{3})$, as in Section 3. Define $p$-adic analogs as follows:

(5.1)

| $\alpha_p$ | $C_{p,\alpha_p}$ | $K_{p,\alpha_p}$ |
|---|---|---|
| 111 | $1/6C_p$ | $(1+p^{-1/3})^3/6K_p$ |
| 21 | $1/2C_p$ | $(1+p^{-1/3})(1+p^{-2/3})/2K_p$ |
| 3 | $1/3C_p$ | $(1+p^{-1})/3K_p$ |
| $1^2 1$ | $1/pC_p$ | $(1+p^{-1/3})^2/pK_p$ |
| $1^3$ | $1/p^2 C_p$ | $(1+p^{-1/3})/p^2 K_p$ |

Here

$$C_p = 1 + p^{-1} + p^{-2} \tag{5.2}$$

$$K_p = \frac{(1 - p^{-5/3})(1 + p^{-1})}{1 - p^{-1/3}}. \tag{5.3}$$

are normalizing factors, chosen so that the column-sums $\sum C_{p,\alpha_p}$ and $\sum K_{p,\alpha_p}$ are both 1. The values of $C_{p,\alpha_p}$ are taken from e.g. [10, page 589] while the $K_{p,\alpha_p}$ are taken from [5, Table 5.1]. Note that both the $C_{p,\alpha_p}$ and the $K_{p,\alpha_p}$ are actually given by procedures uniform in $\alpha_p$. For example, the ratios $(K_{p,\alpha_p} K_p)/(C_{p,\alpha_p} C_p)$ are all given by replacing a factor $f^e$ in $\alpha_p$ by a factor $(1 + p^{f/3})$, e.g. $1^2 1 \mapsto (1 + p^{1/3})(1 + p^{1/3})$.

We expect that Conjecture 3.1 holds verbatim for these more general $\alpha$, with

$$C_\alpha = \prod_{v \in \mathrm{support}(\alpha)} C_{v,\alpha_v} \tag{5.4}$$

$$K_\alpha = \prod_{v \in \mathrm{support}(\alpha)} K_{v,\alpha_v}. \tag{5.5}$$

Our heuristic continues to make sense in this more localized context. In fact the basic references [5] and [6] are in this context.

## 6. Numerical evidence for the refined conjecture

In this section, we work with twenty different $\alpha$. All of them have support $\{\infty, p\}$ and $\alpha_\infty = +$. We let $p$ run over $\{2, 3, 5, 7\}$ and let $\alpha_p$ run over the five-element set $A$. For simplicity, we compare computation with theory only at $x = 10^7$.

We get the numbers $g_\alpha(10^7)$ by renormalizing the first four lines of [10, Table 11]. We get the numbers $f_\alpha(10^7)$ by a simple procedure based on (2.1). Including these extra 501 fields, with total mass $501/3 = 167$, improves the fit slightly, as would be expected from line $j = 7$ of Table 2.

Table 3 gives the comparison. This table illustrates how the slowness of the convergence $h_\alpha(x)/H_\alpha(x)$ to 1 seems to depend systematically on $\alpha_p$, the totally split case $\alpha_p = 111$ being the slowest. On the other hand, $h_\alpha(x)/H_\alpha^*(x)$ is much closer to 1, with no apparent tendency to be above or below.

TABLE 3. Computation and theory: $p$-adic decomposition ($x = 10^7$)

|  | p=2 | | p=3 | | p=5 | | p=7 | |
|---|---|---|---|---|---|---|---|---|
|  | $h_\alpha(x)$ | $h_\alpha(x)$ | $h_\alpha(x)$ | $h_\alpha(x)$ | $h_\alpha(x)$ | $h_\alpha(x)$ | $h_\alpha(x)$ | $h_\alpha(x)$ |
|  | $H_\alpha(x)$ | $H_\alpha^*(x)$ | $H_\alpha(x)$ | $H_\alpha^*(x)$ | $H_\alpha(x)$ | $H_\alpha^*(x)$ | $H_\alpha(x)$ | $H_\alpha^*(x)$ |
| 111 | 0.706 | 0.9999 | 0.721 | 0.9998 | 0.734 | 0.9994 | 0.740 | 0.9994 |
| 21 | 0.851 | 0.9999 | 0.856 | 1.0001 | 0.858 | 0.9996 | 0.858 | 0.9997 |
| 3 | 0.924 | 1.0001 | 0.923 | 0.9999 | 0.920 | 1.0002 | 0.917 | 1.0008 |
| $1^2 1$ | 0.836 | 0.9996 | 0.835 | 0.9997 | 0.832 | 1.0002 | 0.829 | 0.9985 |
| $1^3$ | 0.909 | 1.0003 | 0.903 | 0.9998 | 0.895 | 1.0006 | 0.891 | 1.0025 |

## 7. Three concluding remarks

By using the five element set $A$, we are falling slightly short of the most locally refined conjecture, because there are different ways a prime can be ramified. Thus, for example, if a field $K$ is $1^3$ at $p = 3$, then there are actually nine possibilities for its 3-adic completion $K_3$; both $C_{3,1^3}$ and $K_{3,1^3}$ should be split into nine parts, in the proportions indicated in the $(p, \alpha_p) = (3, 1^3)$ slot on the following chart:

(7.1)

| $\alpha_p$ | $p = 2$ | $p = 3$ | $p \geq 5$ |
|---|---|---|---|
| $1^2 1$ | $\frac{1}{4}\frac{1}{4}\frac{1}{8}\frac{1}{8}\frac{1}{8}\frac{1}{8}$ | $\frac{1}{2}\frac{1}{2}$ | $\frac{1}{2}\frac{1}{2}$ |
| $1^3$ | $1$ | $\frac{1}{3}\frac{1}{3}\frac{1}{9}\frac{1}{81}\frac{1}{81}\frac{1}{81}\frac{1}{81}\frac{1}{81}\frac{1}{81}$ | $1$ or $\frac{1}{3}\frac{1}{3}\frac{1}{3}$ |

Similarly, if a field $K$ is $1^3$ at $p \geq 5$, there are either 1 or 3 possibilities for $K_3$, according to whether $p \equiv 5$, or 1 modulo 6. In general, a $p$-adic algebra $K_p$ of type $\alpha_p$ contributes to $C_{p,\alpha_p}$ and $K_{p,\alpha_p}$ in proportion to $1/(|\operatorname{Aut}(K_p)||D_{K_p}|)$.

Our basic theoretical references [5] and [6] work over arbitrary number fields as well as over finite extensions of $\mathbf{F}_p(t)$ for $p \geq 5$. We expect that the analog of our fully refined Conjecture 3.1 holds in this extra generality too, but at present there is no corresponding published computational data.

The most important remaining problem is to prove our conjecture. The generally pessimistic discussion on [6, page 124] and [2, page 620] suggests to us that the way may be difficult. However, one ingredient of a proof might be the functional equation of $\xi_{\alpha,S}$ with respect to $s \mapsto 1 - s$, studied in [14], [5], and [11]. Another ingredient might be [12, Thm. 3], which concerns growth of arithmetic functions whose associated Dirichlet series satisfy such a functional equation.

Department of Mathematics, Rutgers University, Piscataway, NJ 08854-8019
*E-mail address*: davrobts@math.rutgers.edu